\documentclass[11pt]{article}
\usepackage{amsmath}
\usepackage{amssymb}
\usepackage{amscd}
\usepackage{latexsym}
\usepackage{theorem}
\usepackage{doublespace}
\usepackage{epsfig}
\usepackage{psfrag}
\usepackage{amsfonts}
\usepackage{enumerate}
\usepackage{euscript}

\newtheorem{theorem}{Theorem}[section]
\newtheorem{lemma}[theorem]{Lemma}

\newtheorem{claim}[theorem]{Claim}

{\theorembodyfont{\rmfamily}
\theoremstyle{plain}
\newtheorem{definition}[theorem]{Definition}

\newtheorem{remark}[theorem]{Remark}
\newtheorem{assumption}{Assumption}
}

\numberwithin{equation}{section}

\input xy
 \xyoption{all}

\newcommand{\F}{{\mathbb F}}

\newcommand{\PP}{{\mathbb P}}

\renewcommand{\dim}{\operatorname{dim}}

\newcommand{\C}{{\mathbb{C}}}

\newcommand{\Z}{{\mathbb{Z}}}
\newcommand{\Q}{{\mathbb{Q}}}

\newcommand{\into}{{\hookrightarrow}}

\newcommand{\hs}{\hspace{3pt}}

\newcommand{\half}{\frac{1}{2}}

\newcommand{\iso}{\cong}

\renewcommand{\O}{\mathcal{O}}

\newcommand{\coeff}{R}

\newcommand{\qed}{\hfill \mbox{$\Box$}\medskip\newline}
\newenvironment{proof}{\noindent {\bf Proof:}}{\qed \par}
\newenvironment{proofofclaim}{\noindent {\bf Proof:}}{$\diamondsuit$ \par\medskip}

\newcommand{\algt}{\mathfrak{t}}

\begin{document}

\begin{spacing}{1.1}

\noindent
{\LARGE \bf $T$-equivariant cohomology of cell complexes and 
the case of infinite Grassmannians}
\bigskip\\

\noindent {\bf Megumi Harada }\footnote{{\tt megumi@math.toronto.edu}} \\
Department of Mathematics, University of Toronto,
Toronto, Ontario M5S 3G3 Canada\smallskip \\
{\bf Andr\'e Henriques }\footnote{{\tt andrhenr@math.mit.edu}}\\
Department of Mathematics, Massachusetts Institute of Technology,
Cambridge, MA 02139 \smallskip\\
{\bf Tara S. Holm }\footnote{{\tt tsh@math.berkeley.edu} \\
\newline \mbox{~~~~}{\it MSC 2000 Subject Classification}:
Primary: 55N91  \hspace{0.1in} Secondary: 22E65, 53D20
\newline \mbox{~~~~}{\it Keywords}:
equivariant cohomology, cell complexes, graphs, affine Kac-Moody groups \newline} \\
Department of Mathematics, University of California,
Berkeley, CA 94720

%\vfill\pagebreak

{\small
\begin{quote}
\noindent {\em Abstract.}
In 1998, Goresky, Kottwitz, and MacPherson showed that for
certain spaces $X$ equipped with a torus action, the equivariant
cohomology ring $H^*_T(X)$ can be described by combinatorial data
obtained from its orbit decomposition. Thus, their theory transforms
calculations of the equivariant topology of $X$ to those of the
combinatorics of the orbit decomposition. Since then, many authors
have studied this interplay between topology and combinatorics.  In this
paper, we generalize the theorem of Goresky, Kottwitz, and MacPherson
to the (possibly infinite-dimensional) setting where $X$ is any
equivariant cell complex with only even-dimensional cells and isolated
$T$-fixed points, along with some additional technical hypothesess on the
gluing maps. This generalization includes many new examples
which have not yet been studied by GKM theory, including
homogeneous spaces of a loop group $LG$. 
\end{quote}
}
\bigskip

\section{Introduction and Background}\label{intro}

The main purpose of this paper is to describe the equivariant
cohomology of homogenous spaces of some affine Kac-Moody groups.
Among these examples are the spaces of based loops, $\Omega K$,
considered as a coadjoint orbit of the extended loop group
\(\widehat{LK} \rtimes S^1\). The space $\Omega K$ is a symplectic
Banach manifold, and the maximal torus $T\subseteq \widehat{LK}
\rtimes S^1$ acts on $\Omega K$ in a Hamiltonian fashion.  This
Hamiltonian system exhibits many properties familiar in symplectic
geometry: its moment image is convex \cite{Ati-Pre,Ter}, and its $T$-fixed
points are isolated.  Hence our motivation is to extend, to these
infinite-dimensional examples, results in finite-dimensional
symplectic geometry that compute equivariant cohomology.  Although the
examples that motivate us come from symplectic geometry, our proofs
rely heavily on techniques from algebraic topology.

We now describe the specific symplectic-geometric results that
we will generalize.
Let $X$ be a compact equivariantly formal\footnote{Equivariant
formality is a technical assumption.  It comes for free in all the
examples we consider.}
$T$-space, where $T =
(S^1)^n$ is a finite-dimensional torus. A theorem of Goresky,
Kottwitz, and MacPherson, which we call 
``the GKM theorem'' in honor of its authors, gives a combinatorial
description of the equivariant cohomology ring $H_T^*(X;\F)$,
where $\F$ is a field of characteristic $0$ \cite{GKM}.   The field
coefficients here are crucial. 

We define the
$k$-stratum\footnote{In the symplectic geometry literature, the space
$X^{(k)}$ is usually referred to as the $k$-skeleton.  We will use
this term in the context of cell complexes, so we are introducing the
word stratum to avoid confusion.} $X^{(k)}$ of $X$ to be
\[
X^{(k)} := \{ x \in X \mid \text{dim}(T \cdot x) \leq k\}.
\]
Thus, the 0-stratum \(X^{(0)}\) is just the set of fixed points
$X^T$. This gives the {\em orbit decomposition} of $X$.  In the GKM
theorem, we pay particular attention to the 0-stratum and the
1-stratum, on which additional hypotheses are made.  Note that
\(X^{(l)} \subseteq X^{(k)}\) for \(l \leq k,\) 
so the fixed points $X^T$ are contained in the 1-stratum.  In the
situation that Goresky, Kottwitz and MacPherson consider, the fixed
points are isolated, the equivariant cohomology $H^*_T(X)$ is a free
$H^*_T(pt)$-module, and the kernel of the restriction map $H^*_T(X)
\to H_T^*(X^T)$ is a torsion submodule. Therefore, this
restriction is an injection into the $T$-equivariant cohomology of the
fixed point set $H_T^*(X^T)$.  It is important to note that the
equivariant cohomology of $X^T$, a finite set of isolated points, is
simply the direct product of polynomial rings:
\[
H_T^*(X^T;\F) = \prod_{p \in X^T} H^*_T(pt;\F) \cong 
\prod_{p \in X^T} \F[x_1,\ldots, x_n],
\]
where the degree of $x_i$ is $2$. The $x_i$ are naturally identified
with characters of $T$ and $H_T^*(pt)$ with the symmetric algebra on
the weight lattice $\Lambda$.

The GKM theorem \cite{GKM} now asserts that the image of $H_T^*(X)$ in
$H_T^*(X^T)$ can be described by simple combinatorial data involving
the orbit decomposition of $X$. The hypotheses on $X$ ensure that the
1-stratum consists only of 2-spheres. These spheres are rotated by $T$
with a weight $\alpha\in \algt^*$ (defined up to sign), and have two
fixed points.  They can only intersect at fixed points. Using
this data, we associate to this $T$-space $X$ a graph $\Gamma =
(V,E)$, with vertex set \(V = X^T,\) and edges joining two vertices if
they are the two fixed points on one of the 2-spheres. Moreover, we
associate to each edge $e$ a weight $\alpha_e$ which is precisely the
weight specifying the action of $T$ on the corresponding
2-sphere. Note
that we may think of $\alpha_e$ as a linear polynomial (i.e. degree
$2$ class) in $H_T^*(pt)$. The image of $H_T^*(X)$
depends only on this graph $\Gamma$ and the isotropy data
$\alpha_e$'s.  The GKM theorem says that
$$
H_T^*(X;\F)\iso \left\{ f:V\to H_T^*(pt;\F)\ \left|
\begin{tabular}{l}\ $f(p)-f(q) = \alpha_e \cdot g$\\ for every
$e=(p,q)\in E$\\ and some $g\in
H_T^*(pt;\F)$\end{tabular}\right.\right\}.
$$
In other words, to each vertex $p$ we assign a polynomial $f(p)$.
These polynomials must satisfy some compatibility conditions according
to the edge weights.  Namely, if $e=(p,q)$ is an edge with weight
$\alpha_e$, then $f(p)-f(q)$ must be a multiple of $\alpha_e$.  We now
give a simple example.

\begin{figure}[h]
\begin{center}
\epsfig{figure=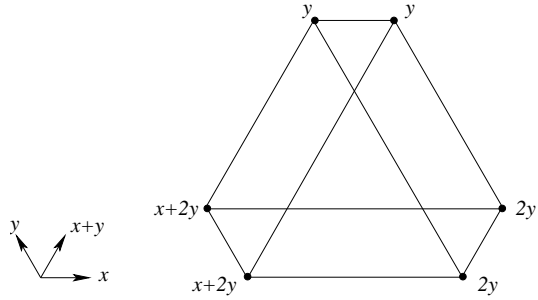,width=2.75in}
\end{center}
\begin{center}
\parbox{5in}{
\caption[coadjoint orbit]{This shows the
GKM graph for a generic coadjoint orbit $\O_\lambda$ of $SU(3)$. The
weights are indicated in the lower left of the figure.  There is a
polynomial attached to each vertex, and the polynomials satisfy the
compatibility conditions, so this does represent an equivariant cohomology
class of $H_T^*(\O_\lambda;\Z)$.}\label{fig:SU3} 
}
\end{center}
\end{figure}

Let $\O_\lambda$ be a generic coadjoint orbit of $SU(3)$.  Then the
maximal torus $T^2$ acts on $\O_\lambda$ by conjugation. There are six
fixed points, the one-stratum is $2$-dimensional, and the associated
GKM graph is shown in Figure~\ref{fig:SU3}.
Other examples of GKM spaces include toric varieties and coadjoint
orbits of any semisimple Lie group.
An identical description can also be given for the equivariant
cohomology of hypertoric varieties \cite{HH04}.
Computations in the equivariant cohomology
ring of flag varieties are closely related to Schubert calculus, and
the GKM description of this ring has added new insights to this
field (see, for example, \cite{Gol,GZ}).

The remainder of the article is organized as follows. In
Sections~\ref{se:injectivity} and \ref{se:GKM}, we generalize the
results of Goresky, Kottwitz and MacPherson to equivariant cell
complexes $X$, possibly infinite-dimensional. We also allow more
general coefficient rings $R$.  There are two results which we must
prove.  The first is the injectivity of $H^*_T(X) \to H^*_T(X^T)$,
which holds when $X$ has only even-dimensional cells; we prove this
in Section~\ref{se:injectivity}.  The second result is
the generalized GKM theorem, combinatorially describing the image of
$H_T^*(X)$ in the cohomology of the fixed points.  This theorem is
more subtle, and requires additional hypotheses on $X$.  We prove this
in Section~\ref{se:GKM}.  In
Section~\ref{generators}, we give a canonical choice of module
generators for $H_T^*(X)$.  In Section~\ref{examples} we discuss the
examples of homogeneous spaces for some affine Kac-Moody groups, which
include the example of $\Omega K$ mentioned above.

\bigskip
\noindent{\bf \em Acknowledgments.}
We are grateful to Allen Knutson for suggesting the problem of looking at the
possible GKM theory for the homogeneous spaces of loop groups and for
teaching the first and second author how to draw GKM pictures. The first
and third authors thank Jonathan Weitsman for many useful discussions.
The first author thanks Robert Wendt for
straightening out all that is twisted (and untwisted) in affine Lie
algebras.

The third author was supported in part by a National Science
Foundation Postdoctoral Fellowship. All authors are grateful for the
hospitality of the Erwin Schr\"odinger Institute in Vienna, where some
of this work was conducted.

\section{The injectivity theorem for cell complexes}\label{se:injectivity}

We show in this section that the equivariant cohomology of $X$ injects
into the equivariant cohomology of its fixed points $X^T$.  Note that
in the category of finite-dimensional symplectic manifolds with
Hamiltonian torus action, this result is a familiar theorem of
Kirwan \cite{Kir}. However, in the more general setting of a cell complex with
$T$-action, we need a separate argument. This result is contained in
Theorem \ref{injectivity}.

We begin with a technical lemma that characterizes the kernel of the
restriction map $H_T^*(Y;R)\to H_T^*(Y^T;R)$.

\begin{lemma}\label{le:borel}
Let $Y$ be a finite dimensional $T$-space with finitely many orbit
types $T/G_i$. Let $R$ be a ring whose torsion is coprime to the
orders of the groups $\pi_0(G_i)$.  Then the kernel of the restriction
map
$$
H_T^*(Y;R)\longrightarrow H_T^*(Y^T;R)
$$
is a torsion $H_T^*(pt;R)$-module.
\end{lemma}

\begin{proof}
Let $Y^T = Y^0\subset Y^1\subset Y^2\subset\cdots\subset Y^m = Y$ be a
filtration of $Y$ so that $Y^i\setminus Y^{i-1}$ has a single orbit
type $T/G_i$. 
We get a spectral sequence from this filtration, with
$$
E_1^{p,q} = H^{p+q}_T(Y^p,Y^{p-1};R) \Longrightarrow H^{p+q}_T(Y;R).
$$
The edge homomorphism $H_T^n(Y;R)\to E_\infty^{0,n} \into E_1^{0,n} =
H_T^n(Y^0;R)$ is the restriction map we are interested in.  We will
prove the lemma by showing that this map becomes an isomorphism after
tensoring with the field of fractions $F$ of $H_T^*(pt;R)$.  This is
because $E_1^{p,q}\otimes F = 0$ for all $p\geq 1$, and so the
spectral sequence collapses after tensoring with $F$.  To show that
$E_1^{p,q} = H^{p+q}_T(Y^p,Y^{p-1};R)$ is torsion for all $p\geq 1$,
we consider the diagram
$$
\xymatrix{
BT \ar[d]_{\pi} & (Y^p\times_T ET, Y^{p-1}\times_T ET) \ar[l] \ar[d] & 
\\
B(T/G_p) &  (Y^p\times_{T/G_p} E(T/G_p) , Y^{p-1}\times_{T/G_p} E(T/G_p))
 \ar[l]_(.75)\varrho \ar[r]^(.68){\cong} & (Y^p/T,Y^{p-1}/T) \\
}.
$$
Let $x\in H^2(B(T/G_p);R)$ be one of the generators.  
Its preimage $\pi^*(x)$ is not a
zero-divisor in $H^2(BT;R)$ by the coprimality assumption.  Since
$Y^p/T$ is finite dimensional, we know that $\varrho^*(x)$ is
nilpotent.  Therefore, $\pi^*(x)$ acts nilpotently on
$$
H^*(Y^p\times_T ET, Y^{p-1}\times_T ET;R) = E_1^{p,q},
$$
completing the proof.
\end{proof}

We now turn our attention to cell complexes.  We say that a space $X$ has
a $T$-invariant cell decomposition if $X$ can be built by
successively attaching cells via $T$-equivariant maps. Each cell has only 
finitely many orbit types. We do not
require the attaching map of a cell to map the boundary to smaller
dimensional cells.  We now state the injectivity result. 

\begin{theorem}\label{injectivity}
Let $X$ be a space with an action of a finite-dimensional torus $T$,
and a $T$-invariant cell decomposition with only even-dimensional
cells, finitely many in each dimension. For any
stabilizer group $G$ of a point, suppose that $R$ is a ring whose
torsion is coprime to the order of the group $\pi_0(G)$.
Let \(\iota: X^T \into X\) denote the inclusion map. Then
the pullback
\[
\iota^*: H^*_T(X;\coeff) \to H^*_T(X^T;\coeff)
\]
is an inclusion.
\end{theorem}

\begin{proof}
We proceed by proving a series of claims.  The main idea of the
proof is to take a non-zero class, restrict it to a {\em finite}
$T$-equivariant cell complex, where we will be able to apply 
Lemma~\ref{le:borel} to conclude injectivity.  Finding the
appropriate finite $T$-equivariant cell complex requires some
knowledge of the module structure of $H_T^*(X;\coeff)$.  We begin
by determining this.

\begin{claim}\label{cl:1}
$H_T^*(X;\coeff)$ is a free $H_T^*(pt;\coeff)$-module, with one
generator\footnote{{\rm If there are infinitely many cells in a given
dimension, all our Theorems still hold, by replacing the phrase ``free
$H^*_T(pt)$-module'' with ``direct product of free rank 1
$H^*_T(pt)$-modules.''}} in degree
$2k$ for each cell of dimension $2k$ for all $k\geq 0$.
\end{claim}

\begin{proofofclaim}
We first show that $H_T^*(X;\coeff)$ is a free
$H_T^*(pt;\coeff)$-module.  Let $X_p$ denote the space built out
of the first $p$ cells.  Consider the cofibration $X_{p}\to
X_{p+1}\to S^{2n}$. By induction, $H_T^*(X_p)$ is evenly graded.
By degree considerations, the long exact sequence of this
cofibration splits into short exact sequences
$$
0\to \tilde{H}_T^*(S^{2n})\to H_T^*(X_{p+1}) \to H_T^*(X_p) \to 0.
$$
Again by induction, $H_T^*(X_p)$ is free, therefore
$H_T^*(X_{p+1}) \cong \tilde{H}_T^*(S^{2n})\oplus H_T^*(X_p)$. Now we
must prove that $\tilde{H}_T^*(S^{2n})$ is a free module with one
generator of degree $2n$.

We will compute $\tilde{H}_T^*(S^{2n};\coeff) =
H^*(S^{2n}\times_T ET,BT;\coeff)$ using the Serre
spectral sequence.  This has
$E_2^{k,\ell}=H^k(BT;\tilde{H}^\ell(S^{2n};\coeff))$, which is
non-zero only when $\ell=2n$.  Thus, this sequence
collapses at the $E_2$ term, since there is only one non-zero row,
and therefore $\tilde{H}_{T}^{*}(S^{2n};\coeff)$ is
$H_T^{*-2n}(pt;\coeff)$.

Now we use the Milnor sequence
\[
0 \to \underleftarrow{\lim}^1\ H^{*-1}(X_i) \to
H^*(\underrightarrow{\lim}\ X_i) \to \underleftarrow{\lim}\ H^*(X_i)
\to 0,
\]
where \(\underleftarrow{\lim}^1\) denotes the first (and only
non-zero) derived functor of the inverse limit functor.  To finish
the argument, we must check that the $\underleftarrow{\lim}^1$ is
zero.  This is true because all the maps $H_T^*(X_{p+1})\to
H_T^*(X_p)$ are surjective, and the inverse limit of surjective
maps is exact. We may now conclude that
$H_T^*(X;\coeff)$ is a free $H_T^*(pt;\coeff)$-module, with one
generator for each cell.
\end{proofofclaim}

Using the above claim, we may think of each
class $\kappa\in H_T^*(X;\coeff)$ as an element
of $H_T^*(pt;\coeff)$,  attached to each cell. Note, however, that the
isomorphism between \(H_T^*(X)\) and the sum \(\bigoplus H_T^*(pt)\) (one
summand for each cell) is not canonical since it relied on choosing
splittings of the projections \(H^*_T(X_{p+1}) \to H^*_T(X_p).\)

Let $\kappa\in H_T^*(X;\coeff)$ be a non-zero equivariant class.
Our goal is to show that $\iota^*(\kappa)\in H_T^*(X^T;\coeff)$ is
also non-zero. As an intermediate step, we restrict our attention
to an appropriate finite $T$-equivariant sub-cell complex $Y$ of
$X$.

\begin{claim}\label{cl:2}
If $\kappa\in H_T^*(X;\coeff)$ is non-zero, then there exists a finite
$T$-equivariant sub-cell complex $Y$, with inclusion $r:Y\into X$, such that
$r^*(\kappa)\in H_T^*(Y;\coeff)$ is non-zero.
\end{claim}

\begin{proofofclaim}
For any $T$-equivariant sub-cell complex $Y$ of $X$, the same argument as
above shows that $H_T^*(Y;\coeff)$ is a free
$H_T^*(pt;\coeff)$-module, generated by its cells. The inclusion
$r:Y\into X$ induces a projection
$H^{*}_T(X;\coeff)\to H^{*}_T(Y;\coeff)$,
which can be interpreted as mapping to zero those generators
corresponding to cells in $X$ but not in $Y$, and is the identity on
the remaining generators.

If $\kappa\in H_T^{*}(X;\coeff;\coeff)$ is a non-zero class, then there
exists a cell $C$ in $X$ such that the equivariant number of $\kappa$
on that cell is non-zero.  Let $Y$ be a finite sub-cell complex
containing that cell.
Then $r^*(\kappa)$ is non-zero in $H_T^*(Y;\coeff)$.
\end{proofofclaim}

Having restricted our attention to a finite $T$-equivariant cell complex
$Y$, we may now apply Lemma~\ref{le:borel}.

\begin{claim}\label{cl:3}
Let $Y$ be a finite $T$-equivariant cell complex with only even-dimensional
cells, and
let \(\iota_Y: Y^T \into Y\) denote the inclusion map. Then
the pullback
\[
\iota_Y^*: H^*_T(Y;\coeff) \to H^*_T(Y^T;\coeff)
\]
is an inclusion.
\end{claim}

\begin{proofofclaim}
By Claim~\ref{cl:1}, $H_T^*(Y;\coeff)$ is a free
$H_T^*(pt;\coeff)$-module.  By Lemma~\ref{le:borel}, the kernel of
$\iota_Y^*$ is a torsion submodule of $H_T^*(Y;\coeff)$.  Therefore,
in this case the kernel must be trivial.
\end{proofofclaim}

We will now show that $\kappa$ must have non-zero image in
$H_T^*(X^T;\coeff)$. Choose $Y$ such that the restriction of $\kappa$ to
$H_T^*(Y;\coeff)$ is non zero. We have the following commutative
diagram
\[
\xymatrix{
H^*_T(X^T;\coeff) \ar[r] & H_T^*(Y^T;\coeff) \\
H^*_T(X;\coeff) \ar[r] \ar[u] & H_T^*(Y;\coeff). \ar[u] \\
}
\]
Since $\kappa$ restricts to be a non-zero
class in $H^*_T(Y^T;\coeff)$ it also restricts to be non-zero in $H^*_T(X^T;\coeff)$.
Hence $H_T^*(X;\coeff)$ injects into $H_T^*(X^T;\coeff)$.
This completes the proof of Theorem~\ref{injectivity}. 
\end{proof}

\begin{remark}
In the above argument, an important step was to show the existence of
the finite-dimensional sub-cell-complex $Y$. It is false, in general,
to claim for a non-zero cohomology class $c$ that there exists a
finite-dimensional sub-cell-complex $Y$ to which $c$ restricts
nontrivially. The problem is that the natural map
\(H^*(\underrightarrow{\lim}\ X_i) \to \underleftarrow{\lim}\ H^*(X_i)\) is
not injective in general. Instead, one has the Milnor sequence \cite{Milnor}
\[
0 \to \underleftarrow{\lim}^1\ H^{*-1}(X_i) \to
H^*(\underrightarrow{\lim} \ X_i)
\to \underleftarrow{\lim} \ H^*(X_i) \to 0.
\]
A simple example where
this $\underleftarrow{\lim}^1$ shows up is in the computation of
\(H^2(K(\Z[\frac{1}{p}],1), \Z) = \Z_p/\Z,\) where the Eilenberg-MacLane
space \(K(\Z[\frac{1}{p}],1)\) is taken to be the direct limit of spaces
homotopy equivalent to $S^1$, and mapping into each other via
cofibrations that induce multiplication by $p$ on $H^1$.

However, in our case, the maps \(X_i \to X_{i+1}\) always induce
surjections in cohomology, and therefore the $\underleftarrow{\lim}^1$ term always
vanishes.
\end{remark}

\section{The GKM theorem for cell complexes}\label{se:GKM}

We now show that the image of the equivariant cohomology of $X$ in
$H_T^*(X^T;\coeff)$ can be identified by simple combinatorial
restrictions involving the $T$-action and the gluing maps. This is the
content of Theorem~\ref{th:GKM}.  The injectivity result of the
previous section is quite general.  We must now make some additional
assumptions on $X$ in order to make this GKM computation.

\begin{assumption}\label{as:1}
The space $X$ can be equipped with a
$T$-invariant cell decomposition, with only even dimensional
cells and only finitely many in each dimension. 
\end{assumption}

\begin{assumption}\label{as:2}
We identify each cell $D^{2n}$ with the unit disc in $\C^n$.
Under this identification, the torus action on $D^{2n}$
is a linear action, given by a group homomorphism \(T \to
T^n,\) where the $T^n$-action on $\C^n$ is the standard action.
\end{assumption}

\begin{assumption}\label{as:3}
The weights \(\{\alpha_i\}\) of the $T$ action on each cell $D^{2n}$ are
pairwise relatively prime as elements of the polynomial ring
$H^*_T(pt;R)\iso R[x_1,\dots,x_k]$, where $k=\dim (T)$. In other
words, if \(\alpha_i|\gamma\) for all $i$, then \(\prod \alpha_i\Big|
\gamma,\) where $\gamma\in H_T^*(pt;R)$.
Moreover, the $\alpha_i$ are not zero divisors. 
\end{assumption}

\begin{assumption}\label{as:4}
Let $W$ denote the cell complex of the first $i-1$ cells, and let
$D^{2n}$ denote the $i^{\mathrm{th}}$ cell.
Let \(\phi: \partial D^{2n} \to W\) be the attaching map for
a cell $D^{2n}$. Then for each $D^2 \subseteq D^{2n}$ corresponding to
an eigenspace of the $T$-action, \(\phi(\partial D^2) \subseteq
W^T,\) i.e. the boundary of each $D^2$ must be mapped to a
fixed point of one of the earlier cells.
\end{assumption}

In Assumption~\ref{as:1}, we could have merely assumed that $X$ was a CW
complex, but this would exclude lots of interesting examples coming
from symplectic geometry, e.g.\ toric varieties. Indeed, in those
examples, the Morse functions used to define the cell structures are
often not Morse-Smale and as a consequence, the cell decompositions
are not CW complexes.

In Assumption~\ref{as:3} we use the identification of the weight
lattice $\Lambda$ with the degree 2 elements \(H^2_T(pt;\coeff).\) Thus it
makes sense to ask that two weights \(\alpha, \alpha' \in \Lambda \cong
H^2_T(pt;\coeff)\) be relatively prime in the ring \(H^*_T(pt;\coeff).\)
Assumptions~\ref{as:2} and \ref{as:3} imply that the $T$ fixed points
of $X$ are isolated, and there is exactly one fixed point for each
cell in the cell decomposition.
By the relative primality in Assumption~\ref{as:3}, we get a
decomposition of each cell $D^{2n}$ into $D^2$'s, corresponding to the
eigenspaces for the $T$-action.

Note that the relative primality assumption also gives a
restriction on the coefficient ring $R$.  Assuming that $\alpha_i$ is
not a zero-divisor in $H_T^*( pt ; R )$ implies that $R$ has torsion
coprime to the orders of the groups $\pi_0(G)$, where $G$ is a
stabilizer group of a point in a cell.  Indeed, if $p$ is a
prime dividing $| \pi_0 ( G ) |$, then $p$ must divide one of the
weights $\alpha_i$ of that cell.  Thus, we may apply the injectivity
result to this $T$-space.

Altogether, these assumptions allow us to define a graph $\Gamma = (V,E)$
associated to $X$. The vertices of $\Gamma$ are the isolated
fixed points of the $T$-action on $X$.  There is an edge connecting
two vertices $p$ and $q$ if these fixed points lie in the
closure of one of the $D^2$'s described in Assumptions~\ref{as:3} and
\ref{as:4}.
Moreover, we associate to each edge the additional datum of a
$T$-weight $\alpha_e$, given by
the weight of $T$ acting on that $D^2$.  We can interpret such a
weight $\alpha_e$ as an element of $H_T^2(pt;\coeff)$.

\begin{definition}
Given a graph $\Gamma=(V,E)$, with each edge $e\in E$ decorated by a
$T$-weight $\alpha_e$, we define the {\em graph
cohomology}\footnote{We apologize for the bad terminology; this is
{\em not} a cohomology theory for graphs.} of $\Gamma$ to
be
$$
H^*(\Gamma) = \left\{ f:V\to H_T^*(pt;\coeff)\ | \ f(p)-f(q)\equiv 0\
({\mathrm{ mod}}\ \alpha_e)\mbox{ for every edge } e = (p,q)\right\}.
$$
\end{definition}

Note that when $\Gamma = (V,E)$ is the graph associated to $X$ as
described above, then this graph cohomology is a subring 
of the equivariant cohomology $H_T^*(X^T;\coeff)$ of the fixed
points $X^T$ of $X$. 
We first prove a Lemma, which computes
the $T$-equivariant cohomology of a 2-sphere. This is the starting
point of the whole discussion.

\begin{lemma}\label{lem:s2coh}
Suppose $T$ acts linearly and non-trivially on $S^2$ with weight
$\alpha$. If $\alpha$ is divisible by $p\in\Z$, assume that $p$ is not
a zero divisor in the coefficient ring $R$. Then the inclusion \(
(S^2)^T=\{N,S\} \into S^2\) induces injections
$
\imath^*: H_T^*(S^2,\{ S\}) \to H_T^*(\{N\})
$
and
$
\jmath^*: H_T^*(S^2) \to H_T^*(\{N,S\}),
$
with images
\[\imath^*(H_T^*(S^2,\{ S\}))=\left\{ g \in H_T^*(\{N\})\ \Big| \ \alpha \mid g
 \right\}\]
and
\[
\jmath^*(H_T^*(S^2)) = \left\{ (f, g) \in H_T^*(\{N\}) \oplus
H_T^*(\{S\})\ \Big| \  \alpha \mid f - g  \right\}.
\]
\end{lemma}

\begin{proof}
The first step is to prove the statement in relative cohomology.
We first consider the case where \(T=S^1,\) acts on
$S^2$ by \(t \cdot z = t^a z\) for some \(a \in \Z,\) not a zero
divisor in $R$.  The cohomology $H^*_T(S^2,\{ S\})$ is equal to
$\tilde{H}^*(S^2\times_{S^1} ES^1/\{ S\}\times_{S^1}ES^1)$, so  we need to
investigate the space $S^2\times_{S^1} ES^1/\{
S\}\times_{S^1}ES^1$. Consider the map
$S^2\times_{S^1} ES^1\to S^2/S^1=[0,1]$.
Its fibers over the endpoints are $BS^1$ and over interior points are
$B (\Z/a\Z)$.  Knowing this, we can write $S^2\times_{S^1} ES^1/\{
S\}\times_{S^1}ES^1$ as $BS^1\cup_f ([0,1]\times B (\Z/a\Z))/(\{
1\}\times B (\Z/a\Z))$ for some map $f:\{0\} \times B (\Z/a\Z)\to BS^1$.  In other
words,  $S^2\times_{S^1} ES^1/\{ S\}\times_{S^1}ES^1$ is the homotopy cofiber
of $f$.  Now consider the long exact sequence of the cofibration
$$
\cdots \to \tilde{H}^*(Cof(f))\stackrel{\imath^*}{\to}
H^*(BS^1)\stackrel{f^*}{\to} H^*(B(\Z/a\Z)) \to
\tilde{H}^{*+1}(Cof(f))\to \cdots.
$$
We know $H^*(BS^1) = R[x]$ with $\deg(x) = 2$.  Because $a$ is not a
zero-divisor in $R$, we also have $H^*(B(\Z/a\Z))) =
R[y]/(ay)$ with $\deg(y) = 2$.  Finally $f^*(x) = y$ up to a unit
since otherwise it would contradict Claim~\ref{cl:1} that
$\tilde{H}^*(Cof(f))$ is evenly graded.  Thus,
$$
    Im(\imath^*) =
    Ker(f^*) =
\left\{ \begin{array}{ll}a\coeff & \mbox{ if $*=2n$ and $n>0$ }\\
0& \mbox{ otherwise}
    \end{array}\right.
$$
which is precisely what we wanted to show.

Now we consider relative case for general $T$. In this case the torus $T$
can always be decomposed as \(T = T' \times S^1,\) where $T'$ acts
trivially on $S^2$ and $S^1$ acts by \(a \in \Z\), as in the
previous case. Then by the K\"{u}nneth theorem, we have the following
diagram

$$
\xymatrix{
H^*_T(S^2,\{ S\}) \ar[r]^{\imath^*} \ar@{=}[d] &
H^*_T(\{N\})\ar@{=}[d] \\
H_{T'}^*(pt) \otimes H^*_{S^1}(S^2,\{ S\})
  \ar[r]^{1\otimes \imath^*_{S^1}} & H_{T'}^*(pt) \otimes
H^*_{S^1}(\{N\}),
}
$$
and so $Im(\imath^*) = H_{T'}^*(pt) \otimes Im(\imath^*_{S^1})$, which
is again what we want.
Now we turn to the non-relative computation.  Consider the following
diagram:
$$
\xymatrix{
0 \ar[r] & H^*_T(S^2, \{S\}) \ar[r] \ar[d]^{\imath^*} & H^*_T(S^2) \ar[r]\ar[d]^{\jmath^*} &
H^*_T(\{S\}) \ar[r]\ar[d] & 0 \\
0 \ar[r] & H^*_T(\{N,S\},\{S\}) \ar[r] & H^*_T(\{N,S\}) \ar[r]
& H_T^*(\{S\}) \ar[r] & 0
}.
$$
Both the top and bottom sequences split, and therefore $Im(\jmath^*) =
Im(\imath^*)\oplus H_T^*(pt)$, where $H_T^*(pt)\to H^*_T(\{N,S\})$ is
the diagonal inclusion.  It is now straight forward to check
\[
\jmath^*(H_T^*(S^2)) = \left\{ (f, g) \in H_T^*(\{N\}) \oplus
H_T^*(\{S\}) \Big| \hs  \alpha \mid f - g  \right\}.
\]
\end{proof}

We need a similar result for a $2n$-sphere.  This is the
technical heart of the proof of Theorem~\ref{th:GKM}.

\begin{lemma}\label{lem:s2ncoh}
Suppose $T$ acts linearly on $S^{2n}$ with $n>1$, and suppose that the weights
$\alpha_1,\dots,\alpha_n$ are pairwise relatively prime over
$H_T^*(pt;R)$. If $\alpha_i$ is divisible by $p\in\Z$,
assume that $p$ is not
a zero divisor in $R$. Then the inclusion \(
(S^{2n})^T=\{N,S\} \into S^{2n}\) induces an injection
$
\jmath^*: H_T^*(S^{2n},\{ S\}) \to H_T^*(\{N\}),
$
with image
\begin{equation}\label{eq:Imj*}
\jmath^*(H_T^*(S^{2n},\{ S\}))=\left\{ g \in H_T^*(\{N\}) \Big|\
\alpha_i \mid g\ \forall i
 \right\}
\end{equation}
\end{lemma}

\begin{proof}
First, we check that the image of $\jmath^*$ is contained in the right
hand side of \eqref{eq:Imj*}.  This is true because we can factor
$\jmath^*$ in $n$ different ways,
$$
H_T^*(S^{2n},\{ S\}) \longrightarrow H_T^*(S^2,\{
S\})\stackrel{\imath^*}{\longrightarrow} H_T^*(\{ N\}).
$$
Thus, by Lemma~\ref{lem:s2coh}, the image of $\jmath^*$ does land in
the right hand side.

Now we show that $\jmath^*$ maps onto the right hand side of
\eqref{eq:Imj*}. We note that $S^{2n}$ is the $n$-fold smash product
$S^{2n} = \bigwedge_{i=1}^n S^2$ of $2$-spheres.  Therefore, it is
possible to use the external cup product to multiply relative
cohomology classes $y_i\in H_T^*(S^2,\{ S\})$ to define a class in
$H_T^*(S^{2n},\{ S\})$.  
Choose a class $g\in
H_T^*(\{ N\})$ satisfying $\alpha_i\mid g$ for all $i$.  By the
relative primality assumption, we conclude that $\prod \alpha_i\mid
g$, and so we can write 
$$
g = \beta \cdot \left(\prod \alpha_i\right). 
$$
Since the $\alpha_i$ are generator of the images
$\imath^*(H_T^*(S^2,\{ S\})) \subseteq H_T^*(\{ N\})$,
by Lemma~\ref{lem:s2coh}, there exist classes $h_i\in
H_T^*(S^2,\{S\})$ 
satisfying $\imath^*(h_i) = \alpha_i$. 
Therefore, the element $\beta \cdot h_1 \smile \ldots \smile
h_n\in H^*(S^{2n},\{ S\})$ satisfies 
$$
\jmath^*(\beta \cdot h_1 \smile \ldots \smile h_n) =
\beta \cdot \imath^*(h_1) \cdots \imath^*(h_n) = g.
$$
Hence $\jmath^*$ is onto the image described in \eqref{eq:Imj*},
completing the proof.
\end{proof}

\begin{theorem}\label{th:GKM}
Let $X$ be a $T$-space satisfying Assumptions~\ref{as:1} through
\ref{as:4}. Then the map
\[
\iota^*: H^*_T(X;\coeff) \to H^*_T(X^T;\coeff)
\]
is an injection, and its image
is equal to the graph cohomology \(H^*(\Gamma) \subseteq
H^*_T(X^T;\coeff),\) i.e.
\begin{equation}\label{eq:GKM}
H^*_T(X;\coeff) \cong H^*(\Gamma).
\end{equation}
\end{theorem}

\begin{proof}
Assumptions~\ref{as:1}-\ref{as:4} imply that the hypotheses of
Theorem~\ref{injectivity} hold, and so we conclude that $\iota^*$ is
an injection.

We now show that the image $\iota^*( H^*_T(X;\coeff) )$ is contained
in $H^*(\Gamma, \alpha)$.  Let $\kappa$ be a class in
$H_T^*(X;\coeff)$, and let \(\iota^*(\kappa)\) be its image in
$H^*_T(X^T;\coeff)$. We denote by \(\iota_p^*(\kappa)\) the
further restriction of $\kappa$ to a single fixed point $p \in
X^T$. To show that \(\iota^*(\kappa)\) is in the graph cohomology, it
suffices to check that for each edge \((p,q) \in E,\) we have the
relation
\[
\iota_p^*(\kappa) - \iota_q^*(\kappa) \equiv 0 \quad (\mathrm{ mod }
\ \alpha_e).
\]
This follows by Lemma~\ref{lem:s2coh} 
from the fact that the restriction of $\kappa$ to the
$S^2$ joining $p$ and $q$ must be an equivariant class in
$H^*_T(S^2;\coeff)$.

We now introduce some notation.  As before, we consider the
filtration of the cell complex $X$ by $X_i$, the set of the first
$i$ cells.    This induces a filtration on the graph
$\Gamma_1\subseteq \Gamma_2\subseteq \cdots \subseteq \Gamma$, 
where $\Gamma_i =
(V_i, E_i)$ has vertices $V_i = (X_{i})^T$ and edges $E_i = \{
(p,q)\in E\ |\ p,q\in V_i\}$.  
We will now prove that the
image of the equivariant cohomology $H_T^*(X)$ is in fact equal to
$H^*(\Gamma)$.  We will prove this by an inductive argument on
the cells. For $X_1 = \{pt\}$, the result is immediate since both
$H_T^*(X_1)$ and $H^*(\Gamma_1)$ are equal to
$H_T^*(pt)$. Now assume the result is known for $X_{i-1}$. We wish
to prove the result for $X_{i}$.

We claim that there is an exact sequence in graph cohomology
\begin{equation}\label{eq:SESgraphcoh}
\xymatrix{
0 \ar[r] & H^*(\Gamma_{i}, \Gamma_{i-1})
 \ar[r]
& H^*(\Gamma_{i}) \ar[r]^{r_i} &
H^*(\Gamma_{i-1}) \ar[r] & 0
}
\end{equation}
where by \(H^*(\Gamma_{i},
\Gamma_{i-1}),\) we mean the {\em relative
graph cohomology} defined by
\[
H^*(\Gamma_{i},\Gamma_{i-1}) :=
 \left\{ f \in H^*(\Gamma_{i})\ |\ f(p) = 0 \mbox{ for
all } p \in V_{i-1} \right\}.
\]
Note that this is not a general property of graph cohomology, but
only holds for those graphs coming from cell complexes. The relative 
graph cohomology consists of exactly those
elements in $H^*(\Gamma_{i})$ whose supports are
concentrated on \(V_{i} \backslash V_{i-1}.\) The map \(r_i:
H^*(\Gamma_{i}) \to H^*(\Gamma_{i-1})\) is
given by the restriction map \(f \mapsto f \vert_{V_{i-1}}.\) The
kernel of $r_i$ is \(H^*(\Gamma_{i},\Gamma_{i-1})\) 
by definition. Therefore,
to show the exactness of the sequence~\eqref{eq:SESgraphcoh}, it
suffices to show that $r_i$ is surjective. To do this, we will use
the following commutative diagram.
\begin{equation}\label{eq:FiveLemmadiag}
\xymatrix{
0 \ar[r] & H^*(\Gamma_i,\Gamma_{i-1}) \ar[r] &
H^*(\Gamma_i) \ar[r]^{r_i} & H^*(\Gamma_{i-1}) &  \\
0 \ar[r] & H_T^*(X_{i}, X_{i-1}) \ar[r] \ar[u] & H_T^*(X_{i})
\ar[r]\ar[u] & H_T^*(X_{i-1}) \ar[r] \ar[u] & 0 \\
}
\end{equation}
The bottom sequence comes from the long exact sequence of relative
cohomology, which automatically splits into short exact sequences,
as before. We know that the right vertical arrow is an isomorphism
by induction. Since the bottom row is exact, a simple diagram
chase implies that the restriction map $r_i$ is surjective.

We will now show that the isomorphism holds at the level of
$X_{i}$ and $\Gamma_i$, i.e.\ that the middle vertical arrow is an
isomorphism. By the Five Lemma, it suffices to show that the left
vertical arrow is an isomorphism.  This is the content of
Lemma~\ref{lem:s2ncoh}.

Finally, we note that
$$
H_T^*(X) = \underleftarrow{\lim}\ H_T^*(X_i) = \underleftarrow{\lim}\
H^*(\Gamma_i) = H^*(\Gamma),
$$
completing the proof.
\end{proof}

\begin{remark}
It is possible to recover the ordinary cohomology $H^*(X)$ from the
$T$-equivariant cohomology by tensoring out the $H^*_T(pt)$. Namely, 
\[
H^*(X;R) = H^*_T(X;R) \otimes_{H^*_T(pt;R)} R. 
\]
Indeed, the Eilenberg-Moore spectral sequence
$Tor_{H^*_T(pt)}(H^*_T(X),R)
\Rightarrow H^*(X)$, coming from the pullback square
\[
\xymatrix{
X \times ET \ar[r] \ar[d] & ET \ar[d] \\
X \times_T ET \ar[r] & BT
}
\]
collapses since $H^*_T(X)$ is a free $H^*_T(pt)$-module.
\end{remark}

\section{Module Generators of $H_T^*(X)$}\label{generators}

We now make the further assumption that $X$ is a CW complex, 
namely that the attaching
maps glue $2n$ cells onto the $2(n-1)$ skeleton. In this section,
we present canonical generators of $H^*_T(X)$ 
as a $H_T^*(pt, \coeff)$-module. We will move freely between thinking
of cohomology classes as either in $H^*_T(X)$ or in 
the graph cohomology
$H^*(\Gamma)$.  

The proofs of Theorems~\ref{injectivity} and \ref{th:GKM} hold
verbatim if we use the filtration by skeleta $X_p$ of the CW
complex.  Assume by induction that we have generators of
$H_T^*(X_{p-1})$. To extend these to $H_T^*(X_{p})$, consider the
short exact sequence
\[
\xymatrix{ 0 \ar[r] & H^{*}_T(X_{p}, X_{p-1}) \ar[r] &
H^{*}(X_{p}) \ar[r] & H^{*}(X_{p-1}) \ar[r] & 0 }.
\]
First, note that \(H^*_T(X_p, X_{p-1}) \cong H^*_T(\bigvee S^p, *)
\cong \bigoplus H^*_T(S^p, *)\) has a canonical choice of
generators. Indeed, each \(H^*_T(S^p, \{S\})\) has a canonical generator,
namely the class whose restriction in \(H^*_T(\{N\})\) is the product
\(\prod_i \alpha_i\) of the weights of the $T$-action on that sphere. 
The generators of $H^{*}(X_{p-1})$ have a unique lift to
$H^{*}(X_{p})$ because $H^{k}_T(X_{p}, X_{p-1})$ is zero for all
$k<p-1$.  These lifts, along with the images of the chosen generators of
$H^{*}_T(X_{p}, X_{p-1})$, form a canonical set of generators of
$H^{*}(X_{p})$.

For each fixed point $v$, let $C_v$ be the corresponding cell, and
$f_v$ be the corresponding generator of $H^*_T(X)$. Let $f_v(w)$
denote the restriction of $f_v$ at the fixed point $w$. It is
straightforward to check that the $\{f_v\}$ satisfy the following
conditions. 

\begin{enumerate}
\item Each $f_v$ is homogeneous of degree \(dim(C_v).\) 
\item If \(\text{dim}(C_w) < \text{dim}(C_v),\) then \(f_v(w) = 0 \in
H_T^*(pt).\)
\item If \(\text{dim}(C_w) = \text{dim}(C_v),\) and \(w \neq v,\) then
\(f_v(w) = 0 \in H_T^*(pt).\)
\item \(f_v(v) = \prod_{i=1}^{\text{dim}(C_v)/2} \alpha_{i} \in
H_T^*(pt),\) 
where the
$\alpha_i$ are the labels of the edges connecting $v$ to
$\Gamma_{\text{dim}(C_v)/2 - 1}$.
\end{enumerate}

These conditions uniquely characterize the $f_v$. Indeed, let $\{f'_v\}$ be
another set of generators satisfying the above conditions. Write them
as \(f'_v = \sum_w b_{vw} f_w.\) By conditions 2 and 3, we have
\(b_{vw} = 0\) whenever \(dim(w) \leq dim(v), w \neq v.\) By condition
4, \(b_{vv} = 1.\) Now, if \(dim(w) > dim(v),\) then \(b_{vw} = 0\) because
otherwise $f_v$ would not be homogeneous. 

\begin{remark}
In the situations where $X$ is a manifold with
a $T$-invariant Morse function $f$ and the cell decomposition is constructed from
the Morse flow with respect to $f$, then the above construction is the same
as the following: given a fixed point $v$, consider the flow-up
manifold $\Sigma_v$ of codimension \(dim(C_v)\). By Poincar\'e
duality, it represents a cohomology class $f_v$ satisfying exactly
these conditions. 
\end{remark}

We illustrate these generators for some examples in the following section.

\section{Grassmannians and flag varieties}\label{examples}

We now turn our attention to the main examples that motivate the
results in this paper. These are the based polynomial loop spaces
$\Omega K$ of a compact simply connected semisimple Lie group $K$, which are
sometimes called the affine Grassmannians. These fall
into the more general category of examples of homogeneous spaces $G/P$
for an arbitrary Kac-Moody group $G$ (defined over $\C$)
with $P$ a parabolic subgroup.  We will phrase the proofs in this
section in a language that makes sense for this more general setting.

We will first consider in Section~\ref{sec:untwLK} the based loop spaces
$\Omega K$. As a homogeneous space, $\Omega K$ has an interpretation
as a coadjoint orbit of $LK$. 
In this setting, the GKM graph can be embedded in $\mathfrak{t}^*$ as
the image of the 1-stratum under a $T$-moment map. The weights
attached to the edges are encoded by their directions. 
Thus, this parallels
the situation for the finite-dimensional coadjoint orbits. 
Throughout this section, we use the coefficient ring
\(R=\Z.\)

\subsection{Based loop spaces $\Omega K$}\label{sec:untwLK}

We first quickly remind the reader of the definitions of
the main characters in this section. 
The loop group $LK$ of $K$
is the set of polynomial loops
\[
LK := \{ \gamma:S^1\to K \},
\]
where the
group structure is given by pointwise multiplication. 
By ``polynomial,'' 
we mean that the loop is the restriction $S^1 = \{z \in \C:
|z|=1\} \to K$ of an algebraic map \(\C^* \to K_{\C}.\) 
The space of {\em based} polynomial loops is defined by
\[
\Omega K = \{ \sigma \in LK |\
\sigma(1) = 1\},
\]
where, by abuse of notation, $1$ is also the identity element in
$K$. It is this space which is GKM space with respect to an
appropriate torus action.

We first observe that $LK$ acts transitively on $\Omega K$ as
follows. For an element \(\gamma \in LK, \sigma \in \Omega K,\) we have
\begin{equation}\label{eq:LSutonPSut}
(\gamma \cdot \sigma)(z) = \gamma(z) \sigma(z) \gamma(1)^{-1}.
\end{equation}
The last correction factor is required to insure that the new loop
\(\gamma \cdot \sigma\) is a {\em based} loop, i.e. that \((\gamma
\cdot \sigma)(1) = 1 \in K.\) This action is clearly transitive, and the
stabilizer of the constant identity loop is $K$.  Hence we may
identify $\Omega K\cong LK/K$.

It is shown in \cite[8.3]{PS} that 
$\Omega K \cong LK/K$ 
is of the form $G/P$ for the
affine group $G = \widehat{LK_{\C}} \rtimes S^1$. Here, $LK_{\C}$ is
the group of algebraic maps \(\C^* \to K_{\C}.\) The
\(\widehat{LK_{\C}}\) is the universal central extension of $LK_{\C}$,
and the $S^1$ acts on $LK_{\C}$ by rotating the loop. 
The parabolic
$P$ is $\widehat{L^+K_{\C}} \rtimes S^1$, where $L^+K_{\C}$ is the
subgroup of $LK_{\C}$ consisting of maps \(\C^* \to K_{\C}\) that
extend to maps \(\C \to K_{\C}.\) The identification is given by
the action of $LK$ on $G/P$
by left multiplication. Then the stabilizer of the identity is $P \cap
LK$. It is the set of polynomial maps \(\C^* \to K_{\C}\) which
extends over $0$ and sends $S^1$ to $K$. A loop \(\gamma\) in $P \cap
LK$ satisfies \(\gamma(z) = \theta(\gamma(1/\bar{z})),\) where
$\theta$ is the Cartan involution on $K_{\C}$. Therefore, since
$\gamma$ extends over zero, by setting \(\gamma(\infty) =
\theta(\gamma(0)),\) it also extends over $\infty$. But then $\gamma$
is an algebraic map from \(\PP^1\) to \(K_{\C},\) and is therefore
constant since $K_{\C}$ is affine. Hence \(P \cap LK = K.\) 

The relevant torus action on $\Omega K$ is given by left
multiplication by the maximal compact torus $T_G$ in $G$. Note, however,
that the center of $G$ acts trivially. Thus we will restrict our
attention to the action of the maximal torus $T_{K_{ad}}$ of $K_{ad} =
K/Z(K)$ and the extra
$S^1$ that rotates the loops. More explicitly, for \(\gamma \in \Omega
K, t \in T_K,\) and \(u \in S^1,\) 
\[
(t,u) \cdot \gamma(z) = t\gamma(uz)\gamma(u)^{-1}t^{-1}.
\]

\subsection{Kac-Moody flag varieties}

We now need to check that this space of based loops $\Omega K = G/P$ 
satisfies
Assumptions~\ref{as:1}-\ref{as:4} that are the hypotheses Theorem \ref{th:GKM}. In fact, the
argument applies to any homogeneous space $G/P$ of a
Kac-Moody group $G$, and $P$ a parabolic, with the action of the
maximal compact torus $T=T_G/Z(G)$ of $G/Z(G)$. It is shown in
\cite{BilDye, KP, KosKum, Mitch}
that $G/P$ admits a CW decomposition 
\[
G/P = \coprod_{[w] \in W_G/W_P} B\tilde{w}P/P.
\]
Here, $W_G$ and $W_P$ are respectively the Weyl groups of $G$ and of
(the semisimple
part of) $P$, and $\tilde{w}$ is a representative of $w$ in
$G$. Each cell has a single $T$-fixed point \(\bar{w} := \tilde{w}P/P.\)
These cells are $T$-invariant
because $T_G$ is a subgroup of $B$, and the center $Z(G)$ acts
trivially. 
To understand the $T$-isotropy weights at each fixed point, we analyze
the tangent space 
\[
T_{\bar{w}}B\bar{w} = T_{\bar{w}}B\tilde{w}P/P =
\mathfrak{b}/\mathfrak{b} \cap \tilde{w}\mathfrak{p}\tilde{w}^{-1}
= \mathfrak{b}/\mathfrak{b} \cap w \cdot \mathfrak{p}.
\]
Therefore, the tangent space decomposes into 1-dimensional pieces, corresponding
to the roots contained in $\mathfrak{b}$ but not in $w \cdot
\mathfrak{p}$. In particular, the weights are all primitive and
distinct. Now pick a root $\alpha$ in $\mathfrak{b}$ but not in $w \cdot
\mathfrak{p}$. Let $e_{\alpha}, e_{-\alpha}$ be the standard root vectors 
for $\alpha, -\alpha$. Let \(SL(2,\C)_{\alpha}\) be the subgroup of $G$ with
Lie algebra spanned by $e_{\alpha}, e_{-\alpha}$, and
$[e_{\alpha},e_{-\alpha}]$ and 
let $B_{\alpha}$ be the Borel of $SL(2,\C)_{\alpha}$ with Lie algebra
spanned by $e_{\alpha}$ and $[e_{\alpha},e_{-\alpha}]$. Let
\(\tilde{r}_{\alpha} := exp(\pi(e_{\alpha} - e_{-\alpha})/2)\)
represent the element $r_{\alpha}$ of the Weyl group which is reflection along
$\alpha$. The $\alpha$-eigenspace in the cell $B\bar{w}$ is
$B_{\alpha}\bar{w} \cong \C$. Its closure is $SL(2,\C)_{\alpha}
\bar{w} \cong \PP^1$, and the point at infinity is given by
$\tilde{r}_{\alpha} wP/P = r_{\alpha} \cdot \bar{w} =
\overline{r_{\alpha}w}$. This is another $T$-fixed point. 
Therefore $G/P$ satisfies all the assumptions of Theorem~\ref{th:GKM}
for \(R = \Z.\)

The GKM graph associated to $G/P$ has vertices $W_G/W_P$, with an edge
connecting $[w]$ and $[r_{\alpha}w]$ for all reflections $r_{\alpha}$
in $W_G$. The weight label on such an edge is $\alpha$. It turns out
that it is possible to embed this GKM graph in $\mathfrak{t}^*$, the
dual of the Lie algebra of $T$, in such a way that the direction of
each edge is given by its label. To produce this embedding, we pick a
point in $\mathfrak{t}_G^*$ whose $W_G$-stabilizer is exactly $W_P$,
take its $W_G$-orbit, and draw an edge connecting any two vertices
related by a reflection in $W_G$. This graph sits in a fixed level of
$\mathfrak{t}_G^*$ (this is only relevant when $G$ is of
affine type) and can therefore be thought of as sitting in
$\mathfrak{t}^*$, where $\mathfrak{t}$ is the Lie algebra of $T=T_G/Z(G)$. 
Since $R = \Z$ and since all weights are primitive,
checking the GKM conditions on the weights amounts to checking that no
two weights from a given vertex are collinear.

\subsection{Moment maps for $\Omega K$}

So far, we
have only considered spaces of polynomial loops in $K$.  However, our
results still apply to other spaces of loops, such as smooth loops,
1/2-Sobolev loops, etc.  Indeed, the polynomial loops are dense in
these other spaces of loops \cite[3.5.3]{PS}, \cite{Mitch}. By Palais'
theorem \cite[Theorem 12]{Pal}, these dense inclusions are weak
homotopy equivalences, and the same holds for the Borel constructions
\(X \times_T ET.\) The statement of Palais' theorem is unfortunately
only stated for open subsets of vector spaces, but can easily be seen
to hold for arbitrary manifolds by a familiar Mayer-Vietoris argument.

\begin{figure}[h]
\begin{center}
\epsfig{figure=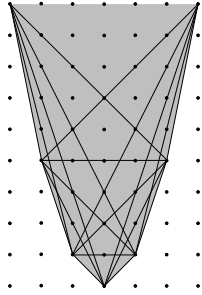,height=1.5in}
\end{center}
\begin{center}
\parbox{5in}{\caption{\small This is the moment
polytope for the $T\times S^1$ action on $\Omega
SU(2)$.}\label{fig:PSut}}
\end{center}
\end{figure}

For $\Omega K$, the embedded GKM graph can be produced as the image of
the 1-stratum under an appropriate $T$-moment map. We first describe
the symplectic structure on $\Omega K$. We write it as a pairing on
$L\mathfrak{k}$. It defines an invariant closed 2-form on $LK$
which descends to $\Omega K$. Let \(X,Y \in L\mathfrak{k}.\) We set
\begin{equation}\label{eq:symplform}
\omega(X,Y) := \int_{S^1} \left<X(t),Y'(t)\right> dt,
\end{equation}
 where $\langle,\rangle$ denotes an invariant bilinear form on $\mathfrak{k}$. The
moment map $\mu: \Omega K \to {\mathfrak t}^*$ 
for the $T$ action is given as follows. Let $X$ denote
an element of $\mathfrak{t}$, and
let \(\gamma \in LK.\) We think of $\gamma$ here as an element of
$LK$, but the formula descends to $\Omega K$. The $X$ component of the
moment map is given by 
\begin{equation}\label{eq:Tmommap}
\mu^X(\gamma) = \int_{S^1} \left<X, \gamma'(t) \gamma(t)^{-1}\right> dt.
\end{equation}
The $S^1$-moment map is given by the {\em energy function},
\begin{equation}\label{eq:S1mommap}
\Phi(\gamma) = \half \int_{S^1} \|\gamma(t)^{-1} \gamma'(t)\|^2 dt.
\end{equation}
The fixed points in $\Omega K$ of the $T\times S^1$-action are exactly
the homomorphisms \(S^1 \to T \subset K,\) and the image of $\Omega K$
under the $T\times S^1$-moment map is the convex hull of the images of
the fixed points \cite{Ati-Pre}. See Figure~\ref{fig:PSut} for the
case \(K=SU(2).\)

\subsection{Loops in $SU(2)$}\label{sec:loops}
%\medskip

We now compute explicitly the ring structure of $H^*_T(\Omega SU(2);
\Z)$ using the moment map graph and the module generators $f_v$ as
constructed in Section~\ref{generators}. In this particular example,
all the restrictions $f_v(w)$ at fixed points $w$ happen to be
elementary tensors in $H^*_T(\{w\})\cong Sym(\Lambda)$, where
$\Lambda=H^2_T(pt)$ is the weight lattice of $T$. This allows us to
use the following convenient notation to represent the classes
$f_v$. On every vertex $w$, we draw a bouquet of arrows $\beta_j\in
\Lambda$ such that $f_v(w)=\prod \beta_j$. The vertices with no arrows
coming out of them carry the class $0$. 

\begin{figure}[h]
\begin{center}
\epsfig{figure=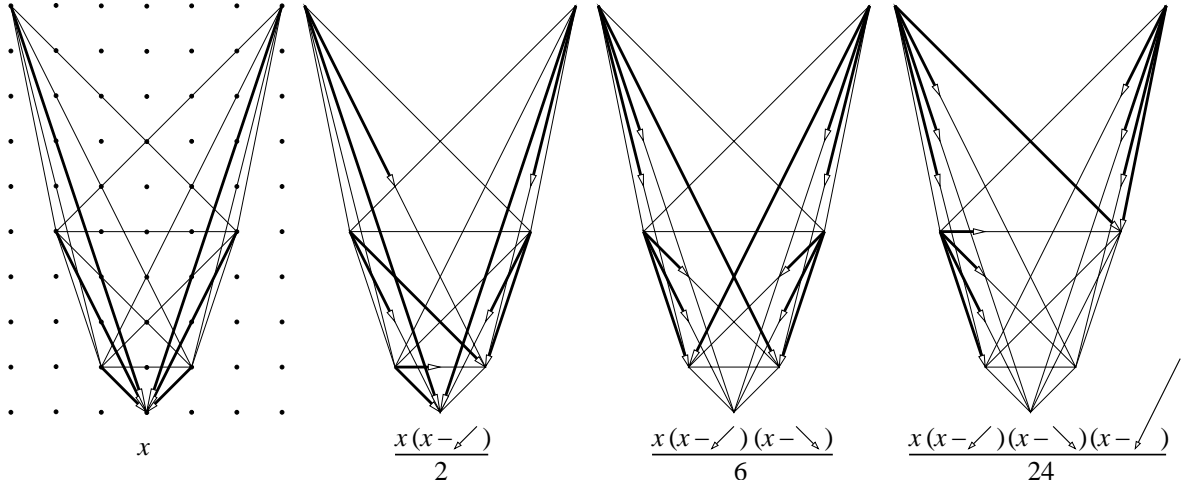,height=2.5in}
\end{center}
\begin{center}
\parbox{5in}{
\caption{
The degree $2$,$4$,$6$, and $8$ generators for $H_T^*(\Omega SU(2); \Z)$. We
draw in the lattice $\Lambda$ in the leftmost
figure.}\label{fig:LGfigures}
}
\end{center}
\end{figure}

The first few module generators are illustrated in
Figure~\ref{fig:LGfigures}. We call $x$ the generator of degree 2, and express
the others in terms of it. The arrows in the expressions denote
elements in $H^2_T(pt)=\Lambda$. 

The map $H_T^*(\Omega SU(2);\Z)\to H^*(\Omega SU(2);\Z)$ is simply the
map that sends the arrows to zero.  And so, by tensoring out the
$H^*_T(pt)$, we recover the well-known fact that 
the ordinary cohomology $H^*(\Omega SU(2); \Z)$ is a divided powers
algebra on a class in degree 2.

If instead, we take the coefficient ring $\Q$, then the cohomology of
$\Omega SU(2)$ is isomorphic to that of $\C P^\infty$.  In fact, there
is a $T^2$ action on $\C P^\infty$ that has the same moment map image
as in Figure~\ref{fig:PSut}.  This $T$-space satisfies 
Assumptions~\ref{as:1} through \ref{as:4} over $\Q$, though not over
$\Z$.  Hence, for this action, $H_T^*(\C P^\infty;\Q)\iso H_T^*(\Omega
SU(2);\Q)$, but this is not true with $\Z$ coefficients.

\subsection{A homogeneous space of type $A_1^{(4)}$}\label{twLoops}

As another example of this type of computation, we let $G$ be the
affine group associated to the Cartan matrix 
\[
\left[
\begin{array}{rr} 2 & -1 \\ -4 & 2 \end{array} \right].
\]
The group is \(\widehat{LSL(3,\C)}^{\Z/2\Z} \rtimes \C^*,\) where the
$\Z/2\Z$-action on $LSL(3,\C)$ is given by precomposition with the
antipodal map \(z \mapsto -z\) on $\C^*$ and composition with the
outer automorphism \(A \mapsto (A^t)^{-1}\) of $SL(3,\C)$.

\begin{figure}[h]
\begin{center}
\epsfig{figure=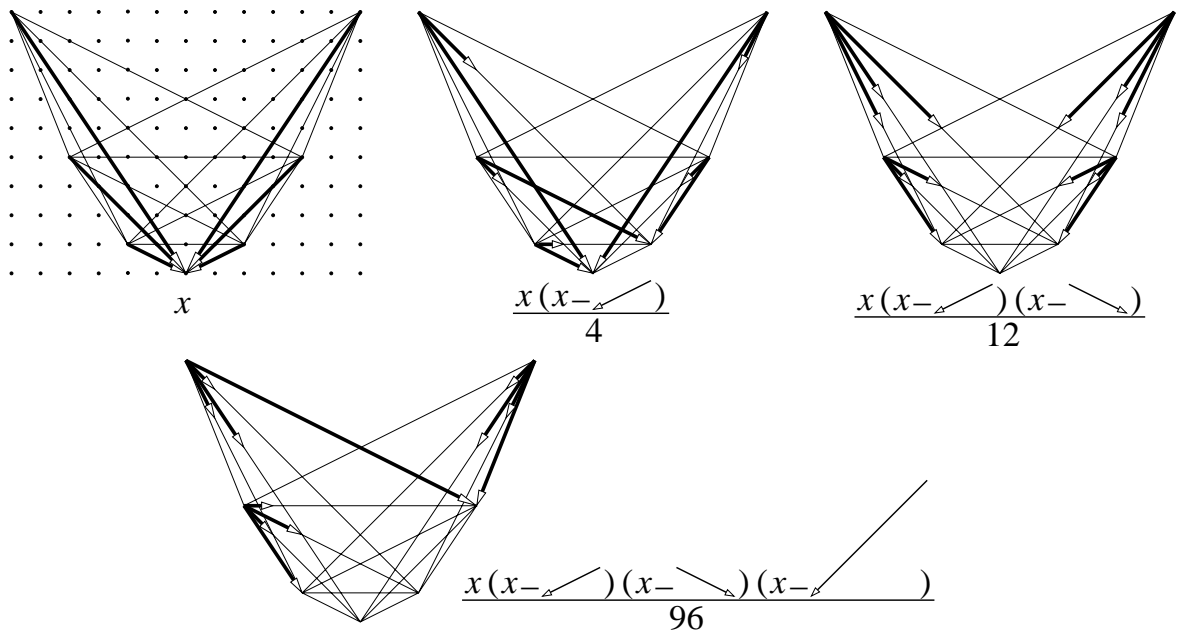,height=3.25in}
\end{center}
\begin{center}
\parbox{5in}{
\caption{The degree 2, 4, 6, and 8 generators
for $H^*_T(G/P; \Z)$.}\label{fig:twLGfigures}
}
\end{center}
\end{figure}

We consider the homogeneous space $G/P$ where the parabolic $P$ has
Lie algebra generated by $\mathfrak{b}$ and the negative of the simple
short root. The degree 2, 4, 6, and 8 module generators in this case
are illustrated in Figure~\ref{fig:twLGfigures}. The denominator in the degree
$n$-th module generator is given by \(n!2^{\lfloor n/2 \rfloor}.\)

\bibliographystyle{plain}
\bibliography{ref}

\end{spacing}

\end{document}